\def\draft{n}
\documentclass{amsart}
\usepackage[headings]{fullpage}
\usepackage{amssymb,epic,eepic,epsfig,amsbsy,amsmath}
\usepackage{amscd,pb-diagram,lamsarrow,pb-lams}

\theoremstyle{plain}

\newtheorem{theorem}{Theorem}\newtheorem{proposition}{Proposition}[section]

\newtheorem{conjecture}{Conjecture}

\theoremstyle{definition}
\newtheorem{definition}[proposition]{Definition}
\newtheorem{problem}{Problem}

\theoremstyle{remark}
\newtheorem{example}[proposition]{Example}
\newtheorem{exercise}[proposition]{Exercise}

\newtheorem{remark}[proposition]{Remark}

\def\printname#1{
        \if\draft y
                \smash{\makebox[0pt]{\hspace{-0.5in}
                        \raisebox{8pt}{\tt\tiny #1}}}
        \fi
}

\newcommand{\psdraw}[2]
         {\begin{array}{c} \hspace{-1.3mm}
        \raisebox{-4pt}{\epsfig{figure=draws/#1.eps,width=#2}}
        \hspace{-1.9mm}\end{array}}

\newlength{\standardunitlength}
\catcode`\@=11
\long\def\@makecaption#1#2{%
     \vskip 10pt

\setbox\@tempboxa\hbox{
       \small\sf{\bfcaptionfont #1. }\ignorespaces #2}%
     \ifdim \wd\@tempboxa >\captionwidth {%
         \rightskip=\@captionmargin\leftskip=\@captionmargin
         \unhbox\@tempboxa\par}%
       \else
         \hbox to\hsize{\hfil\box\@tempboxa\hfil}%
     \fi}
\font\bfcaptionfont=cmssbx10 scaled \magstephalf
\newdimen\@captionmargin\@captionmargin=2\parindent
\newdimen\captionwidth\captionwidth=\hsize
\catcode`\@=12

\def\lbl#1{\label{#1}\printname{#1}}

\def\BN{\mathbb N}
\def\BZ{\mathbb Z}

\def\BQ{\mathbb Q}
\def\BR{\mathbb R}
\def\BC{\mathbb C}

\def\RR{\mathbb R}

\def\calN{\mathcal N}

\def\calL{\mathcal L}

\def\calS{\mathcal S}

\def\a{\alpha}

\def\l{\lambda}

\def\ga{\gamma}

\def\la{\langle}
\def\ra{\rangle}

\def\w{\omega}

\def\e{\epsilon}

\def\d{\delta}

\def\th{\theta}

\def\sub{\subset}

\def\longto{\longrightarrow}
\def\w{\omega}

\def\Om{\Omega}

\def\calB{\mathcal{B}}
\def\bfm{\mathbf{m}}
\def\mul{\mathbf{mul}}
\def\mur{\mathbf{mur}}
\def\mul{\mathbf{mul}}
\def\med{\mathbf{med}}
\def\erfi{\mathrm{Erfi}}
\def\erf{\mathrm{Erf}}

\def\z{\zeta}

\def\gbs{generalized Borel summable}
\def\srb{square root branched}
\def\calK{\mathcal{K}}
\def\calE{\mathcal{E}}

\begin{document}


\title[Resurgence of the Kontsevich-Zagier power series]{
Resurgence of the Kontsevich-Zagier power series}

\author{Ovidiu Costin}
\address{Department of Mathematics \\
         Ohio State University \\ 
         231 W 18th Avenue \\
         Columbus, OH 43210, USA}
\email{costin@math.ohio-state.edu}
\author{Stavros Garoufalidis}
\address{School of Mathematics \\
         Georgia Institute of Technology \\
         Atlanta, GA 30332-0160, USA \\ 
         {\tt http://www.math.gatech} \newline {\tt .edu/$\sim$stavros } }
\email{stavros@math.gatech.edu}

\thanks{The authors were supported in part by NSF. \\
\newline
1991 {\em Mathematics Classification.} Primary 57N10. Secondary 57M25.
}

\date{August 5, 2010 } 


\begin{abstract}
The paper is concerned with the Kontsevich-Zagier formal power series
$$
f(q)=\sum_{n=0}^\infty (1-q)\dots (1-q^n)
$$
and its analytic properties. To begin with, we give an explicit formula
for the Borel transform of the associated formal power series 
$F(x)=e^{-1/(24x)}f(e^{-1/x})$ from which its analytic continuation, its
singularities and their structure can be manifestly determined. This
gives rise to right/left and median summation of the original power series.
These sums, which are well-defined in the open right half-plane are expressed 
by an integral formula involving the Dedekind eta function. The median sum
can also be expressed as a series involving the complex error function.
Moreover, using results of Zagier, we show that the limiting values at
$-1/(2 \pi i \a)$ for rational numbers $\a$ coincide with 
$F(-1/(2 \pi i \a))$.
One motivation for studying the series $f(q)$ is Quantum Topology,
which assigns numerical invariants to knotted 3-dimensional objects.
Our results encourage us to formulate a resurgence conjecture for the
formal power series of knotted objects, which we prove in the case of
the trefoil knot and the Poincare homology sphere, and more generally for 
torus knots and Seifert fibered 3-manifolds. In a subsequent publication we 
will study resurgence for a class of power series that includes the
quantum invariants of the simplest hyperbolic $4_1$ knot.
\end{abstract}

\maketitle

\tableofcontents

\section{Quantum invariants of knotted objects and their puzzles}
\lbl{sec.problem}

\subsection{Introduction}
\lbl{sub.intro2}

The paper is concerned with the Kontsevich-Zagier formal power series
\begin{equation}
\lbl{eq.fq0}
f(q)=\sum_{n=0}^\infty (1-q)\dots (1-q^n)
\end{equation}
and its analytic properties. To begin with, we give an explicit formula
for the Borel transform of the associated formal power series 
$F(x)=e^{-1/(24x)}f(e^{-1/x})$ from which its analytic continuation, its
singularities and their structure can be manifestly determined. This
gives rise to right/left and median summation of the original power series.
These sums, which are well-defined in the open right half-plane are expressed 
by an integral formula involving the Dedekind eta function. The median sum
can also be expressed as a series involving the complex error functions.
Moreover, it is shown using results of Zagier that the limiting values at
$-1/(2 \pi i \a)$ for rational numbers $\a$ coincide with 
$F(-1/(2 \pi i \a))$.
One motivation for studying the series $f(q)$ is Quantum Topology,
which assigns numerical invariants to knotted 3-dimensional objects.
Our results encourage us to formulate a resurgence conjecture for the
formal power series of knotted objects, which we prove in the case of
the trefoil knot and the Poincare homology sphere, and more generally for 
torus knots and Seifert fibered 3-manifolds.
In a subsequent publication we 
will study resurgence for a class of geometrically interesting knotted 
3-dimensional objects that include the simplest hyperbolic $4_1$ knot.

\subsection{Numerical invariants of knotted 3-dimensional objects}
\lbl{sub.intro}

Perturbative quantum field theory assigns numerical invariants (such
as formal power series invariants) to knotted objects. These formal
power series, although they are given by explicit formulas, are typically
factorially divergent, and somehow they are linked to numerical invariants
of knotted objects, such as the Witten-Reshetikhin-Turaev invariants
of 3-manifolds and the Kashaev invariants of knots.

These numerical invariants have poor analytic behavior, satisfy no
known differential equations (linear or not) and the existence of 
asymptotic expansions is a difficult and interesting analytic problem.

In our paper, we formulate a resurgence conjecture for the formal
power series invariants, and show how resurgence solves the numerous
analytic problems, and implies the existence of asymptotic expansions,
and even the presence of exponentially small corrections.

The bulk of our paper consists of a proof of our resurgence conjecture
for the case of the simplest non-trivial
knot, the trefoil ($3_1$), and one of the simplest closed 3-manifolds, 
the Poincar\'e homology sphere. Our results extend without change to torus
knots and Seifert fibered integer homology spheres.

En route, we explain the important notion of resurgence, due to \'Ecalle,
in a self-contained manner.

In a subsequent publication, we will show resurgence of power series
associated to a class of 
geometrically interesting knoted objects, such as the simplest hyperbolic
$4_1$ knot; see \cite{CG1,CG2}. For a detailed discussion of conjectures,
see also \cite{Ga}.

\subsection{TQFT invariants of knotted objects}
\lbl{sub.TQFT}

Let us begin by recalling some of the numerical invariants of
knotted objects. The reader who wishes to focus on the results, may
skip this section, and go directly to Theorem \ref{thm.11} and
Section \ref{sub.results}.

{\em Topological Quantum Field Theory} (TQFT in short) assigns numerical
invariants to knotted 3-dimensional objects. The invariants of 
knots/3-manifolds depend on some additional data, such as a {\em complex
root of unity} $\w$. 
We will denote the numerical invariants by $\phi_{\calK}(\w)$ where 
$\calK$ denotes a {\em knotted object}, that is, a
knot $K$ in 3-space or an integer homology 3-sphere $M$.
In other words, we have a map:
\begin{equation}
\lbl{eq.TQFTphi}
\phi: \text{Knotted objects} \longto \BC^{\Omega}
\end{equation}
where $\Om$ denotes the set of complex roots of unity.
The invariant $\phi_{M}$ is the {\em Witten-Reshetikhin-Turaev invariant} of 
the closed 3-manifold $M$; see \cite{RT,Tu1,Tu2,Wi}.
The invariant $\phi_{K}(e^{2 \pi i/N})$ is the {\em Kashaev invariant} of a 
knot $K$ in 3-space; see \cite{Ka}. Murakami-Murakami showed that 
$\phi_{K}(e^{2 \pi i/N})$ is also equal to the value $\la K \ra_N$ of the 
$N$-th colored Jones polynomial of $K$ (normalized to be $1$ at the unknot), 
evaluated at the $N$-th complex root of unity $e^{2 \pi i/N}$; see \cite{MM}.

The following problem was formulated by Witten (for closed 3-manifolds)
and by Kashaev (for knots).

\begin{problem}
\lbl{prob.1}
Show the existence of asymptotic expansions of the sequence  
$(\phi_{\calK}(e^{2 \pi i/N}))$, and identify the leading terms  
with known geometric invariants; see \cite{Wi,Ka}.
\end{problem}

Unfortunately, the complex-valued function $\phi_{\calK}$,
defined on the set of complex roots of unity, does not seem extend to
a continuous function on the unit circle. Moreover, its asymptotic expansion
around a complex root of unity is unknown, and seems to be a difficult 
analytic problem.

\subsection{Perturbative TQFT invariants of knotted objects}
\lbl{sub.perturbative}

There is an additional formal power series invariant of knotted objects:
\begin{equation}
\lbl{eq.Fpert}
F : \text{Knotted objects} \longto \BQ[[1/x]]
\end{equation}
which is usually thought of as a {\em perturbative expansion} of the quantum
invariants $\phi_{\calK}$. For a homology sphere $M$, $F_M(x)$ is the 
well-known 
Le-Murakami-Ohtsuki invariant (composed with the $\mathfrak{sl}_2$ weight 
system); see for example \cite{LMO} and \cite{Le}.
For a knot $K$, $F_K(x)$ is the Taylor series expansion at $q=e^{-1/x}$ of a 
reformulation of the Kashaev invariant due to Huynh-Le, \cite{HL}. In other 
words, we may write:

\begin{equation}
\lbl{eq.FpertS}
F_{\calK}(x)=\sum_{n=0}^\infty a_{\calK,n} \frac{1}{x^n} \in 
\BQ[[1/x]].
\end{equation}
For every knotted object $\calK$, 
the series $F_{\calK}(x)$ is known to be Gevrey-1 (see \cite{GL})
and in general they are expected to be divergent.

\begin{problem}
\lbl{prob.2}
Show the existence of asymptotic expansions of the sequence
$(a_{\calK,n})$, and identify the leading terms  
with known geometric invariants.
\end{problem}
Thus, we have two types of invariants of a knotted object $\calK$:

\begin{itemize}
\item[(a)] 
the function $\phi_{\calK}: \Om \longto \BC$,
and 
\item[(b)]
the power series $F_{\calK}(x)$. 
\end{itemize}
Using suitable arithmetic completions, in \cite{Ha1,Ha2} Habiro proves that 
either one of the following invariants: 
$F_{\calK}(x)$, $\phi_{\calK}$, the sequence $(\phi_{\calK}(e^{2 \pi i/n}))$,
the sequence $(a_{\calK,n})$, determines the other. We should 
point out that Habiro's proof is in a sense transcendental, of arithmetic
nature. For example, finitely many terms of the sequence $(a_{\calK,n})$ 
{\em cannot} determine $\phi_{\calK}(e^{2 \pi i/3})$.

\begin{problem}
\lbl{prob.3}
Give an analytic proof of Habiro's result.
\end{problem}

Summarizing, we have the following problems:
$$
\divide\dgARROWLENGTH by2
\begin{diagram}
\node[2]{\text{Knotted Objects} }              
\arrow{sw,l}{F}\arrow{se,l}{ \phi }
\\
\node{\text{Formal model} }            
\arrow[2]{e,t,-}{? }\arrow{se,l}{?}
\node[2]{\text{Geometric model} }
\arrow{sw,l}{?}
\\
\node[2]{\text{Asympt.Expansions}}              
\end{diagram}
$$

\subsection{A resurgence conjecture}
\lbl{sub.conjecture}

Despite the apparent analytic difficulties of the series \eqref{eq.fq0}
when $q$ is inside or on or outside the unit circle, and the apparent
factorial divergencies, there seems to be sufficient order
and regularity. Our starting point is the formal power series $F_{\calK}(x)$. 
Let us state the conjecture here, and explain the terms a little later.
For further discussion, see also \cite{Ga}.

\begin{conjecture}
\lbl{conj.1}
For every knotted object $\calK$, 
\begin{itemize}
\item[(a)]
the series $F_{\calK}(x)$ has resurgent Borel transform, 
\item[(b)]
the median sum $S^{\med}_{\calK}$ of $F_{\calK}(x)$ is an analytic function
defined on the right half-plane $\Re(x) >0$, with radial limits
at the points $\frac{1}{2 \pi i}\BQ$ of its natural boundary.
\item[(c)]
Moreover, for $\a \in \frac{1}{2 \pi i}\BQ$, $\a \neq 0$, we have:
\begin{equation}
\lbl{eq.conj1}
S^{\med}_{\calK}\left(-\frac{1}{\a}\right)=\phi_{\calK}(\a).
\end{equation}
\end{itemize}
\end{conjecture}

Our next result shows how resurgence answers the three problems mentioned
above. To state it, recall some standard notation from
asymptotic analysis. For a function $f(x)$ defined a right-hand plane
$\Re(x) >0$, the notation

\begin{equation}
\lbl{eq.Onot}
f(x)=O\left(\frac{1}{x^N}\right)
\end{equation} 
means that there exists positive constants $C$ and $M$
so that $|f(x)| < C/|x|^N$ for all $x$ with $\Re(x) >0$, $|x| > M$.
Furthermore, we say that $f(x)$ is 
{\em asymptotic in the sense of Poincar\'e} to a formal
power series $\hat{f}(x)=\sum_{k=0}^\infty c_k/x^k$ (and write 
$f(x) \sim \hat{f}(x)$) iff for every $N \in \BN$ we have:

\begin{equation}
\lbl{eq.asexp}
f(x)-\sum_{k=0}^{N-1} \frac{c_k}{x^k}= O\left(\frac{1}{x^N}\right).
\end{equation}

\begin{theorem}
\lbl{thm.11}
Assuming Conjecture \ref{conj.1}, it follows that
\begin{itemize}
\item[(a)]
In the interior $\Re(x) >0$, 
\begin{equation}
\lbl{eq.conj2}
S^{\med}_{\calK}(x) \sim F_{\calK}(x) 
\end{equation}
for large $x$.
\item[(b)]
There exist transseries expansions for the sequence $(a_{\calK,n})$ and for
the sequence $S^{\med}_{\calK}(n/(2 \pi i))$.
\item[(c)]
$F_{\calK}(x)$ determines $\phi_{\calK}$ and vice-versa.
\end{itemize}
\end{theorem}
For a definition of a transseries and a proof, see Section \ref{sec.facts}.
Schematically, Conjecture \ref{conj.1} implies the following:

$$
\divide\dgARROWLENGTH by2
\begin{diagram}
\node{\text{Formal model} }            
\arrow{se,r}{\text{Borel transform}}
\node[2]{\text{Geometric model} }
\\
\node[2]{\text{Convolutive model}}            
\arrow{s}\arrow{ne,r}{\text{Laplace transform}}
\\
\node[2]{\text{Transseries Expansions}}              
\arrow{s}
\\
\node[2]{\text{Asympt.Expansions}}              
\end{diagram}
$$

Thus, our Resurgence Conjecture \ref{conj.1} solves at once 
Problems \ref{prob.1}, \ref{prob.2} and \ref{prob.3}
from Sections \ref{sub.TQFT} and \ref{sub.perturbative}.

As a step towards Conjecture \ref{conj.1}, in \cite{GL} Le and 
the first author show that $F_{\calK}(x)$ is a Gevrey power series.

Aside from the applications in Quantum Topology, the conjectured resurgent
series in Conjecture \ref{conj.1} seem to have a different origin than
differential equations. Getting a little ahead of us, the resurgent function
\eqref{eq.Gformula} below does not satisfy any linear (or 
nonlinear) differential equation with polynomial coefficients, as follows
from the structure of its singularities. Resurgence seems to come from
the knotted objects themselves, their combinatorial encodings 
and the exact quantum field theory invariants. This will be investigated 
further in a subsequent publication.

\section{Testing Conjecture \ref{conj.1}}
\lbl{sec.testing}

\subsection{The Zagier-Kontsevich power series}
\lbl{sub.zk}

In the present paper we will
verify the conjecture for the simplest non-trivial knot: the trefoil
$3_1$ (and also for the Poincar\'e Homology sphere; see \ref{sec.appendix}). 
Consider the {\em Kontsevich-Zagier} formal power series

\begin{equation}
\lbl{eq.fq}
f(q) = \sum_{n=0}^\infty (q)_n,
\end{equation}
where the $q$-{\em factorial} $(q)_n$ is defined by
$$
(q)_n=(1-q) \dots (1-q^n)
$$
for $n>0$ with $(q)_0=1$. Although $f(q)$ is {\em not} an analytic function 
of $q$ inside or outside the unit circle, 
it has Taylor series for $q=1$, as well as evaluations at complex roots
of unity. With the notation of Section \ref{sub.perturbative},
we have:
$$
F_{3_1}(x)=e^{-1/x} f(e^{-1/x}).
$$
with $f(q)$ given in \eqref{eq.fq}.
The power series $f(q)$ appears in the beautiful paper of Zagier 
(see \cite{Za}), and was also considered by Kontsevich in a talk at the 
Max-Planck-Institut f\"ur Mathematik in October 1997.
Our basic object of study will be a modified version of $F_{3_1}(x)$, namely,
\begin{equation}
\lbl{eq.Fx}
F(x)=e^{-1/(24 x)} f(e^{-1/x}) \in \BQ[[1/x]]
\end{equation}
with $f(q)$ given in \eqref{eq.fq}.

\subsection{Three models of resurgence in a nutshell}
\lbl{sub.resurgence}

Before we proceed, we need to explain resurgence,
the key aspect of Conjecture \ref{conj.1}.
Resurgence was introduced and studied by \'Ecalle, see \cite{Ec}.
The input of resurgence are formal power series 
and the ouptut are constructible analytic functions in suitable domains, 
which are asymptotic to the original formal power series.
For an extended introduction to resurgence, the reader may also consult 
\cite{DP,D}.

The idea of resurgence is summarized in the following diagram:
$$
\begin{diagram}
\node{F(x) \in \BC[[1/x]] }              
\arrow[2]{e,t,..}{\calS}\arrow{s,l}{ \calB }
\node[2]{\calL^{\bfm}(x)} 
\\
\node{(\calB F)(p), \,\, p=o( 1) }            
\arrow{e,t}{\text{anal. cont.} }
\node{(\calB F)(p) \,\, \text{multivalued} } 
\arrow{e,t}{\text{averaging}}
\node{ (\bfm \calB F)(p) \,\, \text{distribution in} \, \BR^+ }
\arrow{n,r}{ \calL } 
\end{diagram}
$$
and its shorthand version:
$$
\divide\dgARROWLENGTH by2
\begin{diagram}
\node{\text{Formal model}}
\arrow[2]{e,t}{\calS} 
\arrow{se,l}{\calB}
\node[2]{\text{Geometric model}}
\\
\node[2]{\text{Convolutive model}}
\arrow{ne,l}{\calL}
\end{diagram}
$$

Let us explain the terminology here.
\begin{itemize}
\item
The input is a {\em Gevrey}-1 formal power series 
$F(x)=\sum_{n=0}^\infty a_n x^{-n}$. That is, a formal power series
such that there exist constants $C,C'>0$ so that
\begin{eqnarray}
\lbl{eq.gevrey1}
|a_n| \leq C' C^n n!
\end{eqnarray}
for all $n \in \BN$.
\item
The {\em Borel transform} $\calB$ is defined by 
\begin{equation}
\lbl{eq.fborel}
\calB: \BC[[\frac{1}{x}]] \longto \BC[[p]], \qquad
\calB\left(\sum_{n=0}^\infty a_n \frac{1}{x^n}\right)=
\sum_{n=0}^\infty a_{n+1} \frac{ p^n}{n!}
\end{equation}
In other words,
$$
\calB(x^{-n-1})=\frac{p^n}{n!}
$$
If the series $F(x)$ is Gevrey-1, it follows that $(\calB F)(p)$ 
is an analytic function in a neighborhood of $p=0$.
\item
The two horizontal arrows {\em endlessly analytically continue} 
$(\calB F)(p)$ in the complex
plane, minus a discrete set $\calN$ of singularities, as a 
multivalued function. In case the set $\calN$ of singularities of 
$(\calB F)(p)$
is a subset of the real line, one obtains a distribution $(\bfm \calB F)(p)$
on the positive real axis $\BR^+$ by means of an averaging $\bfm$.
This is explained in detail in Section \ref{sub.averaging}. 
\item
The  vertical arrow is the {\em Laplace transform} defined by
$$
\calL \bfm \calB F : \{ x \,\, \Re(x) > c\} \longto \BC,
\qquad (\calL \bfm \calB F)(x) = \int_0^\infty e^{-xp}(\bfm \calB F)(p) dp
$$  
under suitable hypothesis on the growth-rate of $(\bfm \calB F)(p)$ 
for large $p$.
\item
The final horizontal arrow is the {\em generalized Borel transform}
which remembers the constant term of $F(x)$ and is defined by
\begin{equation}
\calS^{\bfm}(F)(x)=a_0+ (\calL \bfm \calB) (F).
\end{equation}
The result is an analytic function defined in a right half-plane.
\end{itemize}

\begin{definition}
\lbl{def.resurgence}
When the above process can be completed, we say that
\begin{itemize}
\item
the formal power series $F(x) \in \BQ[[1/x]]$ is {\em generalized
Borel summable}, (and belongs to the {\em formal model})
\item
its Borel transform $G(p)$ is {\em resurgent}, (and belogs to the
{\em convolutive model})
\item
the resulting function $\calS^{\bfm}(F)(x)$ is {\em analyzable}
(and belongs to the {\em geometric model}).
\end{itemize}
\end{definition}

In what follows, given a \gbs\ series $F(x)$, 
we will denote by $G(p)$ its Borel transform, and by $S^{\bfm}$ its
summation with respect to $\bfm$.

Why is this a reasonable definition? The answer may find an answer in
the following proposition. For a proof, see \cite{Ec} and also the
exposition in \cite{DP,D,CNP1,CNP2,Ml,Ra}.

\begin{proposition}
\lbl{prop.useful}
\rm{(a)} Generalized Borel summation $S^{\bfm}(x)$ coincides with $F(x)$
in case $F(x)$ is analytic in a neighborhood of $x=\infty$:
$$
F(x) = S^{\bfm}(x)
$$
This follows from the following computation
$$
x^{-n-1}=\int_0^\infty e^{-px} \frac{p^n}{n!} dp,
$$
(valid for $x \in \BC$ with $\Re(x) >0$, and  
$n \in \BN$) and the fact that if $F(x)$ is analytic in a neighborhood of 
$x=\infty$, then its Borel transform $G(p)$ is an entire function
of exponential growth, thus the analytic continuation and the averaging
steps do not change $G(p)$, and the Laplace transform reproduces $F(x)$.
\newline
\rm{(b)} If $F(x) \in \BC[[1/x]]$ is \gbs\ with $\bfm$-summation 
$S^{\bfm}(x)$, then for large $\Re(x)$ we have an asymptotic expansion:
$$
S^{\bfm}(x) \sim F(x).
$$ 
\rm{(c)} The set of \gbs\ is an algebra, closed under differentiation 
with respect to $x$.
In particular, if $F(x)$ is a formal solution of a differential 
or difference (linear or not)
equation, then $S^{\bfm}(x)$ is an actual solution
of the equation asymptotic to $F(x)$.
\newline
\rm{(d)}
Generalized Borel summability is a constructive approach, which has 
applications to the numerical approximation of analyzable functions
which are asymptotic to divergent formal power series. 
See for example, the method of {\em truncation to least term}
of factorially divergent series in \cite{CK}.
\end{proposition}
In other words, in analysis we have the following diagram:

$$
\divide\dgARROWLENGTH by2
\begin{diagram}
\node[2]{\text{ODE/PDE} }              
\arrow{sw}\arrow{se}
\\
\node{\text{Formal model} }            
\arrow[2]{e}
\node[2]{\text{Geometric model} }
\arrow{sw}
\\
\node[2]{\text{Convolutive model}}              
\arrow{nw}
\end{diagram}
$$

\subsection{Statement of the results}
\lbl{sub.results}

Let us postpone the remaining definitions to 
Section \ref{sec.laplace}. Our main theorem is the following.

\begin{theorem}
\lbl{thm.1}
\rm{(a)}
The formal power series $F(x)$ of \eqref{eq.Fx} has resurgent Borel transform.
\newline
\rm{(b)}
The median summation
$S^{\med}$ defined on $ \{x \in \BC | \Re(x) >0\}$ 
extends to the points $ \frac{1}{2 \pi i} \BQ$ of its natural boundary 
and for all $\a \in \BQ$, $\a \neq 0$, we have:
\begin{equation}
\lbl{eq.evalcrit}
S^{\med}\left(-\frac{1}{2 \pi i \a}\right)=e^{\pi i \a/12}f(e^{2 \pi i\a}).
\end{equation}
\end{theorem}

The reader may compare Equation \eqref{eq.Fx} that defines the formal
power series $F(x)$ with Equation \eqref{eq.evalcrit} that evaluates the
median summation $S^{\bfm}(x)$ of $F(x)$.

A side bonus is the following precise description of the Borel and Laplace
transforms of $F(x)$. Among other things, it explains why we are using the
median Laplace transform, and identifies the Laplace transforms of our
paper with several functions considered by Zagier in \cite{Za}.

Let $G(p)$ denote the formal Borel transform of the power series $F(x)$
of \eqref{eq.Fx}. Recall the definition of the 
Dedekind eta function $\eta$ and the modified eta 
function $\tilde\eta$ from Section \ref{sec.identities}.
Let $\chi(\cdot)$ be the {\em unique primitive character of conductor} $12$. 
In other words, we have:
\begin{equation}
\lbl{eq.chi}
\chi(n)=
\begin{cases}
1 & \text{if} \,\, n \equiv 1,11 \bmod 12 \\
-1 & \text{if} \,\, n \equiv 5,7 \bmod 12 \\
0 & \text{otherwise.}
\end{cases}
\end{equation}

\begin{theorem}
\lbl{thm.2}
\rm{(a)}
The Borel transform $G(p)$ of $F(x)$ is given
by:
\begin{equation}
\lbl{eq.Gformula}
G(p)=\frac{3\pi}{2\sqrt{2}} 
\sum_{n=1}^\infty \frac{\chi(n) n}{(-p + n^2 \pi^2/6)^{5/2}}.
\end{equation}
$G(p)$ is an analytic double-valued function on $\BC-\mathcal{N}$, 
with singularities in the set $\calN \subset \BR^+$:
\begin{equation}
\lbl{eq.S}
\mathcal{N}= \frac{\pi^2}{6} \Big\{ n^2\,\, | \,\, n \in \BN,
n \equiv 1,5,7,11 \bmod 12 \Big\}.
\end{equation}
\rm{(b)} The left and right summations $S^{\mul}$ and
$S^{\mur}$ are given by:
\begin{eqnarray}
\lbl{eq.lapmul}
S^{\mul}: \{\arg(x) \in (-5 \pi/2,\pi/2) \} \qquad
S^{\mul}(x)= \sqrt{3} x^{3/2} \int_{\gamma_{\e+\arg(x)}}
\frac{\eta(2 \pi i z)}{(x-z)^{3/2}} dz -1 \\
\lbl{eq.lapmur}
S^{\mur}: \{\arg(x) \in (- \pi/2,5 \pi/2) \},
\qquad S^{\mur}(x)= \sqrt{3} x^{3/2} \int_{\gamma_{\e-\arg(x)}}
\frac{\eta(2 \pi i z)}{(x-z)^{3/2}} dz -1 
\end{eqnarray}
where $\gamma_{\theta}$ denotes the ray $\{r e^{i \theta} | r\geq 0 \}$
in the complex plane from $0$ to infinity.
\newline
\rm{(c)}
For every reality-preserving average $\bfm$ (defined in Section 
\ref{sub.averaging}), 
the summation $S^{\bfm}$ is independent of $\bfm$, agrees with the
median summation and is given by:
$$
S^{\med}: \{ x\in \BC| \,\,\, \Re(x) >0 \} \longto \BC
$$
\begin{eqnarray*}
S^{\med}(x) &=& \frac{1}{2} \left( S^{\mur}(x) + S^{\mul}(x) \right)
\\
&=&  \frac{\sqrt{3} x^{3/2}}{2} \left( 
\int_{\gamma_{\e+\arg(x)}}
\frac{\eta(2 \pi i z)}{(x-z)^{3/2}} dz +
\int_{\gamma_{\e-\arg(x)}}
\frac{\eta(2 \pi i z)}{(x-z)^{3/2}} dz \right) - 1. 
\end{eqnarray*}
Moreover,
$$
\overline{S^{\med}(x)}=S^{\med}(\overline{x}).
$$
\rm{(d)} The associated Dirichlet series, defined by 
\begin{equation}
\lbl{eq.dir1}
\d: \{ \Re(x) >0 \} \longto \BC, \qquad
\d(x)=\frac{1}{2} \left( S^{\mur}(x) - S^{\mul}(x) \right)
\end{equation}
equals to:
\begin{equation}
\lbl{eq.dir2}
\d(x) =
i \sqrt{2}(\pi x)^{3/2} 
\tilde\eta (2 \pi i x)
\end{equation}
where $\tilde\eta$ is given in \eqref{eq.etatildef}. 
$\d$ is a lacunary series, with natural boundary $\Re(x)=0$ and with
well-defined radial limits at $\frac{1}{2 \pi i} \BQ$.
\newline
\end{theorem}
The above theorem gives a formula for the median Laplace transform of $G(p)$
in terms of the modified $\eta$-function $\tilde\eta$. Our last theorem
is an alternative formula for the median Laplace transform $S^{\med}$ in terms 
of the {\em complex error function}:

\begin{equation}
\lbl{eq.erfii}
\erfi(x)=\frac{2}{\sqrt\pi} \int_0^{x} e^{t^2} dt
\end{equation}
The complex error function is related to the better known {\em error function}
$$
\erf(x)=\frac{2}{\sqrt{\pi}} \int_0^{x} e^{-t^2} dt
$$
by
$$
\erfi(x)=\frac{\erf(ix)}{i}.
$$
Erfi is an entire odd function of $x$, with asymptotic expansion for large
$x$ with $\arg(x) \in (0,\pi)$ of the form:
$$
\erfi(x) \sim -i + \frac{e^{x^2}}{x \sqrt{\pi}} \sum_{k=0}^\infty
\frac{(2k-1)!!}{(2x^2)^k}
$$
where $(2k-1)!!=1.3\dots(2k-1)$ and $(-1)!!=1$. 
See for example, \cite[Sec.2.2]{O} or \cite[Sec.2]{Lb}.
Consider the modified error function 

\begin{equation}
\lbl{eq.erfiimod}
\calE(x)=e^{-x^2} x^3 \erfi(x) -\frac{x^2}{\sqrt{\pi}}.
\end{equation}
Notice that
\begin{equation}
\lbl{eq.erfbound}
\calE(x)=o(1)
\end{equation}
for large $x$ with $\arg(x) \in (-\pi/4,\pi/4)$.
 

\begin{theorem}
\lbl{thm.lap}
The median Laplace transform of the series $G(p)$ of Equation \eqref{eq.borel}
is given by: 
\begin{eqnarray}
\lbl{eq.laplace2a}
S^{\med}(x)&=&  \frac{12 \sqrt{3}}{\pi^{3/2}}
\sum_{n=1}^\infty \frac{\chi(n)}{n^2}
\calE\left(n \pi \sqrt{\frac{x}{6}}\right) -1.
\end{eqnarray}
\end{theorem}
Notice that Equation \eqref{eq.erfbound} implies that
the series \eqref{eq.laplace2a} is uniformly convergent in
the open right half-plane.

\subsection{Plan of the proof}
\lbl{sub.plan}

Our goal is to motivate, introduce and use resurgence in a relatively
self-contained
fashion. In Section \ref{sec.boreltransform} we compute explicitly the
Borel transform of the series \eqref{eq.Fx} using as input 
the generating function of the Glaisher's numbers, studied by Zagier.
The trigonometric form of this generating function quickly leads, via
a residue computation, to formula \eqref{eq.Gformula} 
for the Borel transform $G(p)$ of Theorem \ref{thm.2}.
This formula is an example of what we call a ``\srb\  function''.
In particular, this implies the existence of the 
analytic continuation of $G(p)$ and locates its singularities.

In Section \ref{sec.laplace} we discuss at length the notion of averaging,
following the work of \'Ecalle, and give several examples of averages. 
Averaging leads to a Laplace transform,
which in general depends on the averaging itself. Our key Proposition
\ref{prop.average} shows that if $G(p)$ is \srb, then all reality-preserving
averages coincide with the median average.
Since our singularities are placed at the positive real numbers, and $G(p)$
is \srb , the difference between the left and right averages is a Dirichlet
series, as we show in Proposition \ref{prop.average}.
We end this section by giving explicit formulas for the median Laplace
transform in terms of the Dedekind $\eta$-function and in terms of the
complex error function, proving Theorem \ref{thm.lap} and part of Theorem
\ref{thm.2}.

In Section \ref{sec.dirichlet} we study the associated Dirichlet series of 
our problem, which turns out to be a modified Dedekind $\tilde\eta$-function.
Zagier's identity and modularity imply the existence of radial
limits of our Dirichlet series. This concludes the proof of
Theorems \ref{thm.1} and \ref{thm.2}.

In Section \ref{sec.facts} we explain how resurgence implies the existence
of asymptotic (and more generally, transseries) expansions of sequences.
In particular, we give a proof of Theorem \ref{thm.11}.

Finally, in Section \ref{sub.extension}, we point out that our results
apply without change to the power series $F_{\calK}(x) \in \BQ[[1/x]]$
where $\calK$ is a $(2,2p)$ torus link or a Seifert fibered rational homology
sphere.

\subsection{Extensions}
\lbl{sub.extension}

For simplicity, we state Theorems \ref{thm.1} and \ref{thm.2}
for the power series $F(x)$ of \eqref{eq.Fx}.

The proof of Theorem \ref{thm.1} works without change 
for the formal power series of $F_{\calK}(x)$ where $\calK$ is a torus link
$(2,2p)$ or a Seifert-fibered homology sphere. In all those cases, 
\begin{itemize}
\item
the Borel transform is \srb , 
\item
the singularities of $G(p)$ are a finite union of sets of the form
$$
\calN=\frac{\pi^2}{\beta} 
\Big\{ n^2 \,\,\, | \,\,\, n \in \BN, \chi(n) \neq 0 \Big\}
$$
for some quadratic character $\chi$. 
\item
the associated Dirichlet series is nearly modular of weight $1/2$,
\item
radial limits of the Dirichlet series exist, 
and Zagier's identity and modularity holds.
\end{itemize}
For an example of the Poincar\'e homology sphere, see the Appendix.
In forthcoming work \cite{CG2} we will prove Conjecture \ref{conj.1} 
for a class of geometrically interesting knotted objects $\calK$ that
include the simplest hyperbolic $4_1$ knot.

\subsection{Acknowledgment}
An early version of this paper was presented by the second author
in a conference in Columbia University, around the Volume Conjecture, 
in the fall of 2005 and spring of 2006.
The second author wishes to thank the organizers of the conference
for their wonderful hospitality.

\section{The Borel transform $G(p)$ of $F(x)$}
\lbl{sec.boreltransform}

In this section we compute the formal Borel transform $G(p)$
of the power series $F(x)$ of \eqref{eq.Fx}. 

\begin{theorem}
\lbl{thm.borel}
We have:
\begin{equation}
\lbl{eq.borel}
G(p)=\frac{3\pi}{2\sqrt{2}} 
\sum_{n=1}^\infty \frac{\chi(n) n}{(-p + n^2 \pi^2/6)^{5/2}}.
\end{equation}
\end{theorem}

Consequently, $G(p)$ is resurgent, and its analytic continuation 
is double-valued in $\BC-\calN$ with singularities in $\calN$, defined as
in Equation \eqref{eq.S}.

\begin{proof}
Let us define a sequence
$(a_n)$ by:
\begin{equation}
\lbl{eq.f}
F(x)=e^{-1/(24x)} f(e^{-1/x})=\sum_{n=0}^\infty  
\frac{a_n}{24^n } \, \frac{1}{x^n}.
\end{equation}
Our sequence $(a_n)$ coincides with Zagier's $(T_n/n!)$ from \cite[Eqn.4]{Za},
where $(T_n)$ are the {\em Glaisher's $T$-numbers}.
In \cite{Za}, Zagier proves that
the {\em Glaisher's $T$-numbers} 
are given by the generating series
\begin{equation}
\lbl{eq.tn}
\sum_{n=0}^\infty \frac{a_n n!}{(2n+1)!} p^{2n+1}=
\frac{\sin 2p}{2 \cos 3p}=
\frac{\sin p}{1-4 \sin^2 p}
\end{equation}
In the following calculations, it will be convenient to let $H(p)$
denote the formal Borel transform of 
\begin{eqnarray*}
F(x/24) &=& e^{-1/x} f(e^{-24/x}) \\
&=& \sum_{n=0}^\infty  \frac{a_n}{x^n}.
\end{eqnarray*}
It is easy to check that
$$
G(p)=\frac{1}{24} H\left(\frac{p}{24} \right).
$$
Thus, it suffices to show that
\begin{equation}
\lbl{eq.borell}
H(p)=1296 \sqrt{3} \pi 
\sum_{n=1}^\infty \frac{\chi(n) n}{(-144 p + n^2 \pi^2)^{5/2}}.
\end{equation}
By the definition of $H(p)$, we have:
\begin{eqnarray*}
H(p) &=& \calB \left( 1 + \sum_{n=0}^\infty a_{n+1} \frac{1}{x^{n+1}} \right) 
\\
&=& \sum_{n=0}^\infty \frac{a_{n+1}}{n!} p^n 
\end{eqnarray*}
On the other hand, Equation \eqref{eq.tn} implies that
$$
p+ \sum_{n=0}^\infty \frac{a_{n+1} (n+1)!}{(2n+3)!} p^{2n+3}=
\frac{\sin p}{1-4 \sin^2 p}
$$
thus
$$
\sum_{n=0}^\infty \frac{a_{n+1}}{n!} \frac{n! (n+1)!}{(2n+3)!} p^{2n}=
\frac{1}{p^3} \left( \frac{\sin p}{1-4 \sin^2 p} -p \right).
$$
Since
$$
\sum_{n=0}^\infty 
\frac{(2n+3)!}{n! (n+1)!} p^{n}=
2 \sum_{n=0}^\infty (2n+3)(2n+1) \binom{2n}{n} p^n=
\frac{6}{(1-4p)^{5/2}}
$$
it follows that
$$
H(p)=(f_1 \circledast  f_2)(p)
$$
where 
\begin{eqnarray*}
f_1(p) &=& \frac{1}{p^{3/2}} \left(
\frac{ \sin( p^{1/2})}{1-4\sin^2 (p^{1/2})}-p^{1/2} \right) \\
f_2(p) &=& 
\frac{6}{(1-4p)^{5/2}} 
\end{eqnarray*}
and $\circledast$ denotes the {\em Hadamard product} 
of two formal power series at $p=0$. The latter is the component-wise
product defined by:
$$
\left(\sum_{n=0}^\infty a_n p^n \right) \circledast
\left(\sum_{n=0}^\infty b_n p^n \right)=
\sum_{n=0}^\infty a_n b_n  p^n. 
$$
It is easy to give a contour integral formula for the Hadamard product:
$$
\lbl{eq.Hg}
H(p)=\frac{1}{2 \pi i} \int_{\gamma} f_1(s) f_2\left(\frac{p}{s}\right) 
\frac{ds}{s}
$$
where $\gamma$ is a small circle around $0$. This  will
give an analytic continuation of the Hadamard product.

Observe that the set $\mathcal{N}_1$ of singularities of $f_1(p)$ is
\begin{equation}
\lbl{eq.Sg1} 
\mathcal{N}_1=\left\{ 
\left(2k+ \frac{1}{6} \right)^2\pi^2,
\left(2k+ \frac{5}{6} \right)^2\pi^2,
\left(2k+ \frac{7}{6} \right)^2\pi^2,
\left(2k+ \frac{11}{6} \right)^2\pi^2 \quad \Big| \quad k \in \BN\right \}
\end{equation}
Now, we enlarge the radius $r$ of the circle $\gamma:=\gamma_r$,
and subtract the residues of the integrand at the singular points,
applying Cauchy's theorem. 

The integrand has single poles at the points $\eta \in \mathcal{N}_1$.
By a straightforward calculation we get that the residue $\psi_{\eta}(p)$
of the integrand at $\eta \in \mathcal{N}_1$ is given by
\begin{equation}
\psi_{\eta}(p)= 1296 \sqrt{3} \pi 
\begin{cases}
\displaystyle 
-\frac{1+12k}{(-144p+(1+12k)^2 \pi^2)^{5/2}}  & \text{if} \,\, 
\eta=(1+12k)^2 \pi^2/6 \\
\displaystyle
+\frac{5+12k}{(-144p+(5+12k)^2 \pi^2)^{5/2}} & \text{if} \,\, 
\eta=(5+12k)^2 \pi^2/6 \\
\displaystyle
+\frac{7+12k}{(-144p+(7+12k)^2 \pi^2)^{5/2}} & \text{if} \,\, 
\eta=(7+12k)^2 \pi^2/6 \\
\displaystyle
-\frac{11+12k}{(-144p+(11+12k)^2 \pi^2)^{5/2}} & \text{if} \,\, 
\eta=(11+12k)^2 \pi^2/6 
\end{cases}
\end{equation}
for $k \geq 0$.
The asymptotic behavior  of $\psi_{\eta}(p)$ for large $\eta$ is 
$O(1/\eta^2)$, 
and thus the sum $\sum_{\eta\in\mathcal{N}_1}\psi_{\eta}(p)$ 
converges. 

$\psi_{\eta}(p)$ is double-valued with only one singularity $\eta/4$ of the 
shape:
$$
\psi_{\eta}(p)=c_{\eta} (\eta/4-p)^{-5/2}.
$$
The definition of $\chi$ and the above computation conclude the proof
of Theorem \ref{thm.borel}.
\end{proof}

\begin{example}
As an numerical check, Theorem \ref{thm.borel} implies that
$$
G(p)=54 \sqrt{3}  \frac{L(4,\chi)}{\pi^4} + O(p)
$$
where 
\begin{equation}
\lbl{eq.Lseries}
L(s,\chi)=\sum_{n=1}^\infty \frac{\chi(n)}{n^s}.
\end{equation}
Since 
\begin{equation}
\lbl{eq.valuesL}
L(2n+2,\chi)= \pi^{2n+2} \frac{(-4)^n}{\sqrt{3} (2n+1)!(n+1)}
\left(B_{2n+2}\left(\frac{1}{12}\right)- 
B_{2n+2}\left(\frac{5}{12}\right)\right)
\end{equation}
(see \cite[Eqn.(6)]{Za})
it follows that 
\begin{equation}
\lbl{eq.Gexample}
G(p)=\frac{23}{24} + O(p).
\end{equation}
On the other hand,
\begin{equation}
\lbl{eq.Fseries}
F(x)=1+\frac{23}{24} \,\, \frac{1}{x} +
\frac{1681}{1152} \,\, \frac{1}{x^2} + 
\frac{257543}{82944} \,\, \frac{1}{x^3} + 
\frac{67637281}{7962624} \,\, \frac{1}{x^4} +
 O\left(\frac{1}{x^5}\right), 
\end{equation}
which confirms that the constant term of the Borel transform of $F(x)$
is given by $23/24$, in accordance with \eqref{eq.Gexample}.
\end{example}

\begin{exercise}
\lbl{exer.1}
Using Theorem \ref{thm.borel} and the special values of the L-series
given in \eqref{eq.valuesL}, show that
$$
G(p)=\frac{23}{24} + \frac{1681}{1152} p + \frac{257543}{165888}p^2+
 \frac{67637281}{47775744}p^3 + O(p^4)
$$
in confirmation with the Borel transform of \eqref{eq.Fseries}.
\end{exercise}

\section{The Laplace transform of $G(p)$}
\lbl{sec.laplace}

In this section we compute the Laplace transform $\calL(x)$ of $G(p)$.

\subsection{Analytic continuation, averaging and Laplace
transform}

Given a resurgent function $G(p)$ with singularities in $\BN^+ \subset \BR^+$,
there are three ways to average and take the Laplace transform.
\begin{itemize}
\item
The first way is to use a uniformizing average $\bfm$ of \'Ecalle
in order to get a single valued function $\bfm G$ 
on $\BR^+$. Unfortunately,
this function is not integrable since $\int_0^1 dp/(p-1)^{-5/2}$ does not
exist. So, \'Ecalle applies an acceleration operator to $\bfm G$
and then takes the usual Laplace transform. 
\item
Alternatively, \'Ecalle applies a uniformizing average to the Laplace transform
of the analytic continuation of $G(p)$ along paths that avoid the 
singularities.
The key property is that the set of such paths form a Riemann surface.
\item
The first author converts $G(p)$ to a step-distribution on 
$\BR^+\setminus\BN^+$
and then applies an extended Borel transform $\calB_{\a}$, followed by
an extended Laplace transform. See \cite[Sec.1.3]{C}.
\end{itemize}

Now, we arrive at a subtle point: there are many well-behaved
uniformizing averages. In fact for every probability distribution 
$f \in L^1(\BR)$ with $\int_0^\infty |f(x)| dx=1$. \'Ecalle-Menous construct
a uniformizing average $\bfm_f$; see \cite{EM}.

On the other hand, the first author extended Borel transforms $\calB_{\a}$
are parametrized by $\a \in 1/2+ i\BR$. Of all those Borel transforms the most
useful one is the balanced one $\calB_{1/2}$, which satisfies the key property
of approximation by summation to least term; see \cite{CK}.

In case a formal power series satisfies a generic differential equation 
(linear or not),
all averages $\bfm_f$ agree with the first author's 
balanced $\calB_{1/2}$, as shown in \cite{C}. This is also a consequence
of \'Ecalle's {\em bridge equation}; see \cite{Ec,DP}.

In our case, the series $F(x)$ does not satisfy 
a differential equation. Nevertheless, Proposition \ref{prop.average} shows 
a universality, i.e., independence of averaging.
Before we state the proposition, let us explain what averaging means.

\subsection{What is an averaging?}
\lbl{sub.averaging}

Averages were introduced and studied extensively by \'Ecalle.
Following \'Ecalle-Menous (see \cite{EM}), 
let us consider a multivalued function $G(p)$ defined on $\BC-\calN$,
with at most exponential growth at infinity, and with singularities
on a discrete set $\calN=\{\eta_k \,\, | \,\, k \in \BN\} \sub \BR^+$, where
$\eta_k < \eta_l $ for $k < l$. 

Let us define the relative spacing $\w_k$ $k \in \BN$
of the singularities by $\w_k=\eta_{k+1}-\eta_k$. Thus, we have the picture:

$$
\psdraw{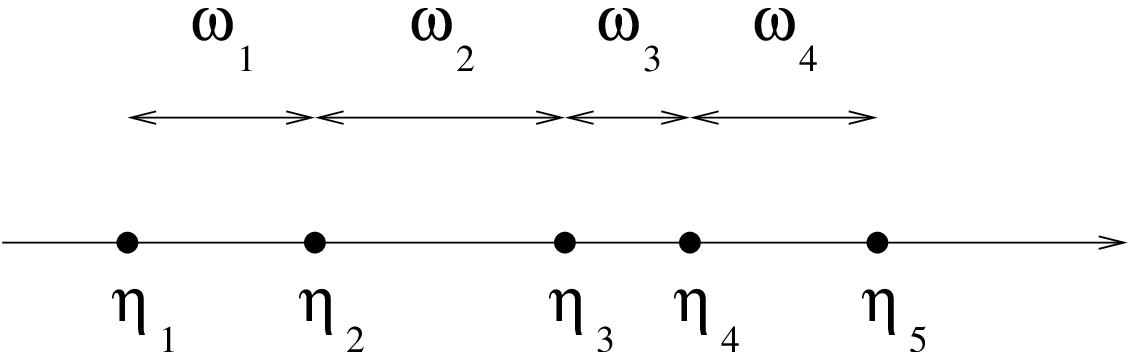}{3in}
$$

With respect to the terminology of \'Ecalle (cf. e.g. \cite{EM}) we have
that $G\in{\textsc{Ramif}}(\RR^+)$ with ramification points $\mathcal{N}$. 
An averaging $\bfm$ is a linear map
$$
\bfm: {\textsc{Ramif}}(\RR^+) \longto 
{\textsc{Unif}}(\RR^+) 
$$
that maps multivalued functions with singularities at $\calN$ to
single-valued distributions on $\BR^+$. Averaging maps depend on a set
of averaging weights.
An {\em averaging weight}
$\bfm^{\varpi}$ is a collection
$$
\{ \bfm^{\varpi}=\bfm^{\varpi_1,\dots,\varpi_r} \quad
| \quad r \in \BN, \varpi_i=\left( \frac{\e_i}{\w_i} \right),
\quad \e_i= \pm, \w_i=\eta_{j+1}-\eta_j \}.
$$
The tuple
$(\epsilon_1,\dots,\epsilon_r)$ is called an {\em address}.
We always assume that for all $r$, we have:
\begin{equation}
\lbl{eq.consistency}
\bfm^{\varpi_1,\dots,\varpi_r}=
\bfm^{\varpi_1,\dots,\varpi_r, \left( \frac{+}{\w_{r+1}} \right)} +
\bfm^{\varpi_1,\dots,\varpi_r, \left( \frac{-}{\w_{r+1}} \right)}
\end{equation}
Recall that $G(p)$ is a multivalued function. For fixed $r \in \BN$ and 
$\varpi_1, \dots, \varpi_r$, we now define a multivalued function
$G^{\varpi_1, \dots, \varpi_r}$ as follows. Let
 
\begin{equation}
\lbl{eq.Gvarpi}
G^{\varpi_1, \dots, \varpi_r}: (\eta_r,\eta_{r+1}) \longto \BC
\end{equation}
denote the analytic continuation of $G$ along a path avoiding the 
singularity at $\eta_i$ from above if $\epsilon_i=+$ and from below otherwise. 
Then, $\bfm G$ is defined by:
$$
\bfm G(p)=\sum_{\e_1=\pm, \dots, \e_r=\pm}
\bfm^{\varpi_1,\dots,\varpi_r} G^{\varpi_1, \dots, \varpi_r}(p),
\qquad
p \in (\eta_r, \eta_{r+1}).
$$

There are several natural properties that are often required for averages
$\bfm$.
Three important properties are:

\begin{itemize}
\item[(P1)]
$\bfm$ preserves reality and has real-valued weights,
\item[(P2)]
$\bfm$ preserves convolution,
\item[(P3)]
$\bfm$ preserves lateral growth.
\end{itemize}

P1 is useful when the input is a power series with real coefficients and
the needed output is an analytic function on the right half-plane $\Re(x) >0$
which takes real values for $ >0$.

P2 is needed for commutation of generalized Borel summability with
multiplication of power series.

P3 is necessary to be able to define Laplace transforms.

Let us give three rather trivial, but useful averages from \cite[p.85]{EM}:
\begin{eqnarray*}
\mathbf{mur}^{\varpi_1,\dots,\varpi_r} =
\begin{cases}
 1 & \text{if} \,\, \e_1=\e_2=\dots=+ \\
0 & \text{otherwise},
\end{cases}
\end{eqnarray*}
\begin{eqnarray*}
\mathbf{mul}^{\varpi_1,\dots,\varpi_r} =
\begin{cases}
 1 & \text{if} \,\, \e_1=\e_2=\dots=- \\
0 & \text{otherwise},
\end{cases}
\end{eqnarray*}
\begin{eqnarray*}
\mathbf{med}^{\varpi_1,\dots,\varpi_r} = 
\begin{cases}
 1/2 & \text{if} \,\, \e_1=\e_2=\dots=+ \\
 -1/2 & \text{if} \,\, \e_1=\e_2=\dots=- \\
0 & \text{otherwise},
\end{cases}
\end{eqnarray*}

$\mul,\mul, \med$ and $\calB_{1/2}$ satisfy P3. $\calB_{1/2}$ and 
$\med$ satisfies P1.

\subsection{The Laplace transform of an averaged function}
\lbl{sub.laplaceave}

Recall that the Laplace transform of a function 
$G(p) \in L^1(\BR^+, e^{-\nu x}dx)$ (for $\nu >0$) with at most
exponential growth at infinity is defined by:
\begin{equation}
\lbl{eq.deflaplace}
\calL G : \{ x \in \BC | \Re(x)>1/\nu \} \longto \BC, \qquad 
(\calL G)(x)= \int_0^\infty e^{-px} G(p) dp.
\end{equation}

If $G(p)$ is defined in a sectorial neighborhood of $0$ (i.e, 
in a set $\{p \in \BC | \arg(p) \in (-\e,\e')$) and is of exponential
growth at infinity, then by moving the 
integration contour it follows that $\calL(x)$ is defined in an enlarged
neighborhood $\{x \in \BC| |x|>1/\nu, \arg(x) \in (-\e'-\pi/2,\e+\pi/2)\}$.

The definition of the Laplace transfrom makes sense in case $G(p)$ is a 
{\em distribution} (e.g.
$G(p)=1/(p-1)$), as was discussed by the first author in \cite[Sec.2]{C}. 
Likewise, we may define the Laplace transform of $\bfm G(p)$:
\begin{equation}
\lbl{eq.laplace3}
(\calL^{\bfm} G)(x)=\int_0^\infty e^{-xp} \bfm G(p)dp
\end{equation}
It turns out that $(\calL^{\bfm}G)(x)$ is an average of line integrals
of $G(p)$ along paths in $\BC-\calN$ that start at $0$ and end at $\infty$.
For example, it is easy to see that
\begin{eqnarray*}
(\calL^{\mur} G)(x) &=& \int_{\gamma_r} e^{-px} G(p) dp \\
(\calL^{\mul} G)(x) &=& \int_{\gamma_l} e^{-px} G(p) dp \\
(\calL^{\med} G)(x) &=& \frac{1}{2} \left( \int_{\gamma_r} e^{-px} G(p) dp
+ \int_{\gamma_l} e^{-px} G(p) dp
\right)
\end{eqnarray*}
where $\gamma_r$ (resp. $\gamma_l$) is a path in $\BC-\mathcal{N}$ from $0$ to 
$\infty$ that turns right (resp. left) at each singularity in $\mathcal{N}$:

$$
\psdraw{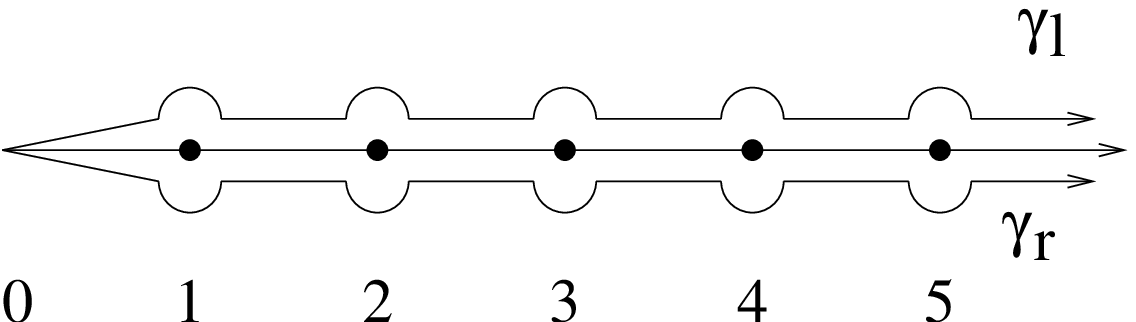}{2in}
$$

Another useful average is $\calB_{1/2}$ of the first author; see 
\cite[Eqn(1.20)]{C}.

In case the multivalued function $G(p)$ is the Borel transform of a formal
power series solution $F(x)$ of a generic
differential (or difference) equation, then 
it is known that the Laplace transforms $\calL^{\bfm}(x)$ for all averages
that satisfy P1, P2, P3 agree. In our case, $F(x)$ is not expected
to satisfy a differential equation (linear or not) with polynomial
coefficients, because the position of singularities (which is an analytic
invariant) is qualitatively different from solutions to differential
equations with polynomial coefficients.
What is a natural average
to consider? The next lemma states that for the singularities of $G(p)$
in Equation \eqref{eq.borel}, the Laplace transform is independent of
the averaging.

To state the lemma, we need some notation. Motivated by Equation 
\eqref{eq.borel}, let us introduce the following definition.

\begin{definition}
\lbl{def.clean}
We will call a multivalued function $G(p)$ {\em \srb } if it is given by
a (Mittag-Leffleg like) absolutely convergent sum:
\begin{equation}
\lbl{eq.clean}
G(p)=\sum_{\eta \in \calN} G_{\eta}(p),
\end{equation}
where $\calN$ is a discrete subset of $\BR^+$,
$$
G_{\eta}(p)= c_{\eta} (\eta-p)^{-k_{\eta}/2},
$$
and $k_{\eta} \in \BN^+$. 
Thus, the support of $c$, $\{\eta \in \calN \,\, | c_{\eta} \neq 0\}$
is the set of singularities of $G(p)$.
\newline
A \srb\  function $G(p)$ has {\em weight} $k$ when $k_{\eta}=k \in \BN^+$.
\end{definition}

\begin{proposition}
\lbl{prop.average}
\rm{(a)} If $G(p)$ is \srb\ of odd weight $k$ and 
$\bfm$ is any \'Ecalle average that preserves P1
(and may or may not preserve P2 or P3),  
then 
\begin{equation}
\lbl{eq.med3}
\bfm G(p)=\sum_{\eta \in \calN} \bfm G_{\eta}(p) 
\end{equation}
where 
\begin{equation}
\lbl{eq.med2}
\bfm
G_{\eta}(p)=\begin{cases}  G_{\eta}(p) & \text{if} \,\, p < \eta \\
0 & \text{if} \,\, p \geq \eta
\end{cases}
\end{equation}
does not depend on $\bfm$. 
\newline
\rm{(b)} For $\Re(x)>0$, we have: 
\begin{equation}
\lbl{eq.laplacelr}
(\calL^{\med} G)(x)=\frac{1}{2} \left( (\calL^{\mur} G)(x) 
+ (\calL^{\mul} G)(x) \right).
\end{equation}
where
$ (\calL^{\mul} G)$ and $(\calL^{\mur} G)$ are defined for $x \in \BC^*$ with
$\arg(x) \in (-5\pi/2,\pi/2)$ and $\arg(x) \in (-\pi/2,5\pi/2)$
respectively.
\newline
\rm{(c)}
In their common domain $\arg(x) \in (-\pi/2,\pi/2)$,
the associated Dirichlet series, is defined by:
\begin{equation}
\lbl{eq.dirichlet2}
 \d(x)=\frac{1}{2} \left( (\calL^{\mur} G)(x) - (\calL^{\mul} G)(x) \right)
\end{equation}
\rm{(d)}
Consequently, we have:
\begin{eqnarray}
\lbl{eq.laplace2c}
(\calL^{\med} G)(x)&=& (\calL^{\mul} G)(x) + \d(x) \\ 
\lbl{eq.laplace2d}
      &=& (\calL^{\mur} G)(x) - \d(x)
\end{eqnarray}
\rm{(e)}
If $c_{\eta} \in \BR$ for all $\eta$, then 
$$
\overline{(\calL^{\mul} G)(x)}=(\calL^{\mur} G)(\overline{x}),
\qquad 
\overline{(\calL^{\mur} G)(x)}=(\calL^{\mul} G)(\overline{x}),
\qquad
\overline{(\calL^{\med} G)(x)}=(\calL^{\med} G)(\overline{x}).
$$
\rm{(f)}
When $k$ is odd, the Dirichlet series is given
by:
\begin{equation}
\lbl{eq.dirichlet}
\d(x)=i \frac{2^{(k-1)/2} \sqrt{\pi} x^{k/2-1}}{(k-2)!!}
\sum_{\eta \in \calN} c_{\eta} e^{-\eta x}  
\end{equation}
where for an natural number $n \in \BN$ we denote 
$(2n+1)!!=1.3.5.\dots (2n+1)$.
\end{proposition}

\begin{proof}
It suffices to consider the case 
$$
G(p)=c_{\eta} (\eta-p)^{-k_{\eta}/2} 
$$
where $k_{\eta}$ is positive integer. Let us fix an average
$\bfm$ of \cite{EM} which is symmetric (i.e., satisfies P2 of Section 
\ref{sub.averaging}) and let us suppose that $\eta=\eta_r$ for some 
$r \in \BN$. Observe that $G(p)$ is not singular for $p \in [0, \eta_r)$.
Equation \eqref{eq.consistency} implies that 
$$
G^{\varpi_1, \dots, \varpi_s}(p)=G(p), \qquad
p \in  (\eta_s,\eta_{s+1})
$$
for $s < r$. On the other hand, for $p \in (\eta_r,\eta_{r+1})$,
the two analytic continuations of the square root differ only in sign;
thus,
$$
G^{\varpi_1,\dots,\varpi_r, \left( \frac{+}{\w_{r+1}} \right)}(p)=
-G^{\varpi_1,\dots,\varpi_r, \left( \frac{-}{\w_{r+1}} \right)}(p),
$$
which together with the symmetry condition P2 imply that
$\bfm G(p)=0$ for $p \in (\eta_r,\eta_{r+1})$, and in fact for
$p > \eta_r=\eta$. This proves (a).

Part (b) follows from Section \ref{sub.laplaceave}.

Parts (d), (e) follow from (b) and (c).

The definition of the Dirichlet series implies that
$$
\d(x)= c_{\eta} \int_{C_{\eta}} (\eta-p)^{-k_{\eta}/2} e^{-px} dp
$$
where $C_{\eta}$ is a loop (Hankel contour) from $+\infty$, $\arg(p)=0$
to $+\infty$, $\arg(p)=2\pi$ which goes once around $\eta$, 
oriented counterclockwise. A residue calculation implies (f).
\end{proof}

\subsection{A formula for the Laplace transform of $G(p)$}
\lbl{sub.formulalaplace}

In this section we will prove of part (c) of Theorem \ref{thm.2} and
Theorem \ref{thm.lap}.
We will use the Dedekind $\eta$ function as in \eqref{eq.etaf}.
Recall the contours $\gamma_\th$ from Theorem \ref{thm.2}.
Recall also that $S^{\mul}(x)$ denotes the Laplace transform of $\mul G(p)$.

\begin{theorem}
\lbl{thm.laplace}
For $x \in \BC$, $x \neq 0$, $\arg(x) \in (-5\pi/2,\pi/2)$, we have:
\begin{equation}
\lbl{eq.laplace}
S^{\mul}(x) 
= \sqrt{3} x^{3/2} \int_{\gamma_{\e+\arg(x)}} \eta(2 \pi i z)
\frac{dz}{(x-z)^{3/2}} - 1.
\end{equation}
\end{theorem}

This proves Equation \eqref{eq.lapmul} of Theorem \ref{thm.2}. 
\eqref{eq.lapmur} is completely analogous.

\begin{proof}
We have:
\begin{xalignat*}{2}
S^{\mul}(x) &= \int_{\gamma_l} e^{-px} G(p) dp  & \\
&= \frac{3 \pi}{2 \sqrt{2}}  \sum_{n=1}^\infty \chi(n) n 
\int_{\gamma_l} \frac{e^{-px}}{(-p+n^2 \pi^2/6)^{5/2}}  dp & 
\text{by Thm \ref{thm.borel}} \\
&=  \frac{3 \pi}{2 \sqrt{2}} \sum_{n=1}^\infty \chi(n) n 
\left( \int_{\gamma_l}  \frac{2 x}{3}
\frac{e^{-px}}{(-p+n^2 \pi^2/6)^{3/2}}  dp - 
\frac{2 \cdot 6^{3/2}}{3 n^3 \pi^3} \right)
& \text{by integration by parts} \\
&
= \sqrt{3} x \sum_{n=1}^\infty \chi(n) 
\int_{\gamma_l} \frac{e^{-n^2 \pi^2 qx/6}}{(-q+1)^{3/2}}  dq
+ C
& \text{by a change of variables}\,\, 6p=n^2 \pi^2 q 
\end{xalignat*}
where
\begin{equation*}
C=:
-6 \sqrt{3} \frac{1}{\pi^2} \sum_{n=1}^\infty \frac{\chi(n)}{n^2} 
=-6 \sqrt{3} \frac{1}{\pi^2} L(2,\chi) = -1
\end{equation*}
where the last Equality follows from Equation \eqref{eq.valuesL}.
Since
\begin{xalignat*}{2}
\sqrt{3} x \sum_{n=1}^\infty \chi(n) 
\int_{\gamma_l} \frac{e^{-n^2 \pi^2 qx/6}}{(-q+1)^{3/2}}  dq
&= 
\sqrt{3} \sum_{n=1}^\infty \chi(n)
\int_{\gamma_{\e+ \arg(x)}}
\frac{e^{-n^2 \pi^2 z/6}}{(-z/x+1)^{3/2}}  dz
& \text{by a change of variables} \,\, z=qx \\
&=  \sqrt{3} x^{3/2} \int_{\gamma_{\e+\arg(x)}} \eta(2 \pi i z)
\frac{dz}{(x-z)^{3/2}} 
& \text{by} \,\, \eqref{eq.etaf}
\end{xalignat*}
the result follows.
\end{proof}

We now give a proof of Theorem \ref{thm.lap}. 

\begin{proof}(of Theorem \ref{thm.lap})
Recall the complex error function $\erfi(x)$ and its modification $\calE(x)$
from Equations \eqref{eq.erfii} and \eqref{eq.erfiimod}.
The median Laplace transform is the average of the left and right Laplace
transform. Moreover, a calculation shows that for $x>0$ we have:

\begin{equation}
\lbl{eq.erfi2}
\int_{\gamma_l} \frac{e^{-xp}}{(1-p)^{5/2}}dp=
-\frac{2}{3} -\frac{4x}{3} -\frac{4}{3} i \sqrt{\pi} e^{-x} x^{3/2}
+ \frac{4}{3} \sqrt{\pi} e^{-x} x^{3/2} \erfi(\sqrt{x}).
\end{equation}
Replacing $\gamma_l$ by $\gamma_r$ has the effect of replacing $i$ by $-i$
in the above equation. Thus, the median integral, which also coincides
with the {\em principal value integral}, is given by:

\begin{equation}
\lbl{eq.erfi3}
\int_{\gamma_m} \frac{e^{-xp}}{(1-p)^{5/2}}dp=
-\frac{2}{3} -\frac{4x}{3} 
+ \frac{4}{3} \sqrt{\pi} e^{-x} x^{3/2} \erfi(\sqrt{x})=
-\frac{2}{3} + \frac{4}{3} \sqrt{\pi} \calE(\sqrt{x})
\end{equation}
where $\gamma_m=1/2(\gamma_l+\gamma_r)$. On the other hand, the proof
of Theorem \ref{thm.laplace} implies that

$$
S^{\med}(x)=
\sqrt{3} x \sum_{n=1}^\infty \chi(n) 
\int_{\gamma_m} \frac{e^{-n^2 \pi^2 qx/6}}{(1-q)^{3/2}}  dq -1.
$$
Using Equation \eqref{eq.erfi3}, the result follows.
\end{proof}

\section{A Dirichlet series $\d(x)$ associated to $F(x)$}
\lbl{sec.dirichlet}

\subsection{A formula for a Dirichlet series $\d(x)$ associated to $F(x)$}
\lbl{sub.dir}

In this section we identify the associated Dirichlet series $\d(x)$
of the \gbs\ power series $F(x)$ of \eqref{eq.Fx} with the 
Eichler integral $\tilde\eta$ of the Dedekind $\eta$-function
given by \eqref{eq.etatildef}. In particular,
using Zagier's identity (see \eqref{eq.strange}) and a modular property
of one of Zagier's functions (see \eqref{eq.evalg}), allows us to prove
the existence of radial limits at complex roots of unity and to finish
the proof of Theorems \ref{thm.1} and \ref{thm.2}.

\begin{proposition}
\lbl{prop.dir}
\rm{(a)}
The Dirichlet series associated to $F(x)$ is given by:
\begin{equation}
\lbl{eq.laplace2b}
\d: \{x \in \BC | \Re(x) >0\} \longto \BC, \qquad
\d(x) =
i \sqrt{2}(\pi x)^{3/2} 
\tilde\eta (2 \pi i x)
\end{equation}
\rm{(b)} $\d$ is a lacunary series with natural boundary the line $\Re(x)=0$.
\newline
\rm{(c)} $\d$ has radial limits at $\frac{1}{2 \pi i}\BQ$ given by:
\begin{equation}
\lbl{eq.evaldir}
\d\left(-\frac{1}{2 \pi i \a} \right)=
\z_{24}^3 \a^{-3/2} \phi(-1/\a)
\end{equation}
for all $\a \in \BQ$, $\a \neq 0$, where $\phi$ is a function of Zagier 
from \eqref{eq.phi} and $\z_k=e^{2 \pi i/k}$.
\end{proposition}

\begin{proof}
Theorem \ref{thm.borel} gives that
$$
G(p)=\frac{3 \pi}{2  \sqrt{2}} 
\sum_{n=1}^\infty \frac{\chi(n) n}{(- p + n^2 \pi^2/6)^{5/2}}.
$$
Part (f) of Proposition \ref{prop.average} implies that the associated
Dirichlet series is given by:
\begin{eqnarray*}
\d(x) &=& 
i \sqrt{2}(\pi x)^{3/2} \sum_{n=1}^\infty \chi(n) n e^{-\pi^2 n^2 x/6} \\
&=& 
i \sqrt{2}(\pi x)^{3/2} 
\tilde\eta (2 \pi i x)
\end{eqnarray*}
where the last equality follows from the definition of $\tilde\eta$ in
\eqref{eq.etatildef}. This proves (a).

(b) follows from \cite{Ma}. In other words, $\d(x)$ cannot be analytically
continued beyond the line $\Re(x)=0$. In general, lacunary series need
not have radial limits at points of their natural boundary. Our series,
however, has radial limits at rational multiples of $1/(2 \pi i)$.

(c) follows from Equation \eqref{eq.laplace2b} and Zagier's
identity \eqref{eq.strange} below.
\end{proof}

\subsection{Proof of of Theorem \ref{thm.1}}
\lbl{sub.proofthm1}

We are finally in a position to finish the proof of Theorem \ref{thm.1}. 

\begin{theorem}
\lbl{thm.eval}
With the notation as in Theorem \ref{thm.1}, 
for all $\a \in \BQ$, $\a \neq 0$, we have:
$$
S^{\med}\left(-\frac{1}{2 \pi i \a} \right)=\phi \left( \a \right).
$$
\end{theorem}

\begin{proof}
Equations \eqref{eq.lapmul}, \eqref{eq.lapmur} and the definition
of Zagier's $g$-function of Equation \eqref{eq.g}
imply that for $\a \in \BQ-\{0\}$ we have:
\begin{equation}
\lbl{eq.calLg}
g(\a)=
\begin{cases} 
S^{\mul}(-1/(2 \pi i\a)) & \text{if} \,\, \a >0 \\
S^{\mur}(-1/(2 \pi i\a)) & \text{if} \,\, \a < 0. \\
\end{cases}
\end{equation}

Let us assume $\a \in \BQ$, $\a >0$ (the other case is analogous).
We have:
\begin{xalignat*}{2}
S^{\med}\left(-\frac{1}{2 \pi i \a} \right) 
&= S^{\mul}\left(-\frac{1}{2 \pi i \a} \right)
+ \d\left(-\frac{1}{2 \pi i \a} \right)
  & \text{by \eqref{eq.laplace2d}}\\
&= S^{\mul}\left(-\frac{1}{2 \pi i \a} \right) +\z_{24}^3 \a^{-3/2}
\phi\left(-\frac{1}{\a} \right) 
  & \text{by Prop. \ref{prop.dir} (c)} \\
&= g(\a) -(i \a)^{-3/2} \phi(-1/\a) & \text{by \eqref{eq.calLg}} \\
&= \phi(\a) & \text{by \eqref{eq.evalg}}
\end{xalignat*}
\end{proof}

\section{Identities from Zagier's paper}
\lbl{sec.identities}

In this section we collect several definitions, notations and results from 
Zagier's paper \cite{Za}, for the convenience of the reader. 
Zagier defines a function
\begin{equation}
\lbl{eq.phi}
\phi: \BQ \longto \BC, \qquad \phi(\a)=e^{\pi i \a/12} f(e^{2 \pi i \a})
\end{equation}
which evaluates at complex roots of unity the series 
$f(q)$ of \eqref{eq.fq}.
Zagier considers the following formal power series in $\BZ[[q]]$:
\begin{eqnarray*}
(q)_\infty &=& \prod_{n=1}^\infty (1-q^n) \\
&=& \sum_{n=-\infty}^\infty (-1)^n q^{n(3n+1)/2} \\
&=& \sum_{n=1}^\infty \chi(n) q^{(n^2-1)/24} \\
H(q) &=& \sum_{n=1}^\infty \chi(n) n q^{(n^2-1)/24}
\end{eqnarray*}
as well as the corresponding analytic functions for $q=e^{2 \pi i z}$,
$\Im (z) >0$:

\begin{eqnarray}
\lbl{eq.etaf}
\eta(z) &=& e^{\pi i z/12} \prod_{n=1}^\infty (1-e^{2 \pi i n z})=
\sum_{n=1}^\infty \chi(n) e^{\pi i n^2z /12} \\
\lbl{eq.etatildef}
\tilde\eta(z)&=& 
\sum_{n=1}^\infty \chi(n) n e^{\pi i n^2z /12}
\end{eqnarray}

$\eta(z)$ is the famous Dedekind $\eta$ function, a modular form of weight 1/2,
and $\tilde\eta(z)$ is an 
{\em Eichler integral} of the Dedekind $\eta$ function.
Although $\tilde\eta$ is not a modular form,  Zagier proves 
that $\tilde\eta$ has radial limits to the rational points 
$z \in \BQ \subset \BR$ of its natural boundary.

Zagier's identity (coined ``the strange identity'' by Zagier himself) 
\cite[Eqn.7]{Za} identifies 
the radial limits of $\tilde\eta$ with $\phi$ for $\a \in \BQ$:
\begin{equation}
\lbl{eq.strange}
\phi(\a)=-\frac{1}{2} \tilde\eta(\a).
\end{equation}

At the last two pages of his seminal paper, 
Zagier introduces a $C^\infty$ function

\begin{equation}
\lbl{eq.g}
g: \BR \longto \BC, \qquad g(x)=\int_0^\infty (z-x)^{-3/2} \eta(z) dz,
\end{equation}
where $\eta(z)$ is the {\em Dedekind} $\eta$ function defined by
\eqref{eq.etaf}. Zagier states that
$g(x)$ is real analytic everywhere except at $x=0$ and whose derivatives 
at $0$ are given by
$$
g^{(n)}(0)=(-\pi i/12)^n n! a_n,
$$
where
\begin{equation}
\lbl{eq.f2}
F(x)=e^{-1/(24x)} f(e^{-1/x})=\sum_{n=0}^\infty  
\frac{a_n}{24^n } \, \frac{1}{x^n}.
\end{equation}
Moreover, for $\a \in \BQ$, we have:
\begin{eqnarray}
\lbl{eq.invg}
g(\a)&=&(i \a)^{-3/2} g(-1/\a) \\
\lbl{eq.evalg}
\phi(\a)+(i\a)^{-3/2} \phi(-1/\a)&=& g(\a) \,\, \text{for} \,\, a \in \BQ 
\end{eqnarray} 

In other words, for $h \to 0$ we have:
\begin{eqnarray}
\notag
g(h) & \sim & \sum_{n=0}^\infty \frac{g^{(n)}(0)}{n!} h^n \\
\notag
& = & \sum_{n=0}^\infty \left( -\frac{ \pi i}{12} \right)^n a_n h^n \\
\lbl{eq.gg3}
&=& e^{ 2 \pi i h/24} f(e^{2 \pi i h})
\end{eqnarray}
where the last equality follows from Equation \eqref{eq.f}.  
In \cite[Eqn.6]{Za} Zagier gives a closed formula for the Taylor 
coefficients $(a_n/24^n)$ of $F(x)$:

\begin{eqnarray*}
\frac{a_n}{24^n} &=& 6  \frac{(-6)^n}{(n+1)!} 
\left(B_{2n+2}\left(\frac{1}{12}\right)- 
B_{2n+2}\left(\frac{5}{12}\right)\right) \\
&=&
\frac{1}{2 \sqrt{3} (\pi/6)^2 (2 \pi^2/3)^n} \frac{(2n+1)!}{n!}
L(2n+2,\chi) 
\end{eqnarray*}
where 
\begin{equation}
\lbl{eq.Lseries2}
L(s,\chi)=\sum_{n=1}^\infty \frac{\chi(n)}{n^s}.
\end{equation}
Consider now the Borel transform
$$
G(p)=\sum_{n=0}^\infty \frac{a_{n+1}}{24^{n+1} n!} p^n
$$
of $F(x)$. To simplify notation, let us write

\begin{equation}
\lbl{eq.betas}
G(p)=\sum_{n=0}^\infty b_n p^n
\end{equation}
instead. Then, we have:

\begin{eqnarray}
\lbl{eq.Gformula1}
b_n &=& 6  \frac{(-6)^{n+1}}{(n+2)!n!} 
\left(B_{2n+4}\left(\frac{1}{12}\right)- 
B_{2n+4}\left(\frac{5}{12}\right)\right) \\
&=&
\lbl{eq.Gformula2}
\frac{4 \pi^2}{\sqrt{3} (2 \pi^2/3)^n} \frac{(2n+3)!}{(n+1)!n!}
L(2n+4,\chi) 
\end{eqnarray}

Since 
$$
\frac{(2n+3)!}{(n+1)!n!} \sim 4^n n^{3/2}
\left(\ga_0 + \frac{\ga_1}{n} + \frac{\ga_2}{n^2}
+
\dots \right)
$$
for computable constants $\ga_j$, 
and since $L(2n+4)=1+O(5^{-2n})$ for every $M$,
it follows that the coefficients of the Borel transform have an asymptotic
expansion of the form:

\begin{equation}
\lbl{eq.Gformula3}
b_n \sim \left(\frac{6}{\pi^2}\right)^n n^{3/2}
\left( c_{1,0} + \frac{c_{1,1}}{n} + \frac{c_{1,2}}{n^2} +\dots \right)
\end{equation}
for computable constants $c_{1,l}$ for $l \in \BN^+$.
Disassembling the L-series into its monomial parts, 
Equations \eqref{eq.Gformula2}  reveals a {\em transseries expression} for
the coefficients of the Borel transform:

\begin{equation}
\lbl{eq.Gformula4}
b_n \sim
\left(\frac{6}{\pi^2}\right)^n n^{3/2}
\sum_{l=0}^\infty\sum_{k=1}^\infty \frac{c_{k,l}}{n^l k^n}
\end{equation}
for a doubly indexed series of resurgence monomials $(6/\pi^2)^n n^l k^n$,
and for computable constants $c_{k,l}$. Notice that the resurgence monomials
form a well-order set of order type $\w^2$.

\section{Resurgence implies transseries expansions}
\lbl{sec.facts}

Let us examine more carefully the asymptotic equations from the last section.
Although Equation \eqref{eq.Gformula2} makes sense,
the asymptotic series in \eqref{eq.Gformula3} is factorially divergent. 
In view of this, one cannot naively make sense of Equation \eqref{eq.Gformula4}
since for example $1/2^n$ is a monomial which is (exponentially)
smaller than any of the monomials $1/n^l$ for all $l$. In order to reach the
monomial $1/2^n$ we would have to subtract the infinite series of all previous
monomials $1/n^l$ for $l \in \w$, and this series is factorially divergent.
What we need is a way to subtract the whole series at once.
It is at this point that resurgence is needed to make sense of the formal
series in \eqref{eq.Gformula}. 

Recall that the singularities of
$G(p)$ are included in the set $ \lambda \BN^+$ where $\lambda=\pi^2/6$.

Fix a small positive angle $\theta$ and for every $k \in \BN^+$
draw the rays $L_k=k\l \e^{i\theta}\BR^+$ 
from $k \lambda$ to infinity along the direction of $\theta$.
Assume that $\theta$ is chosen so that the rays $L_k$ are distinct:
$$
\psdraw{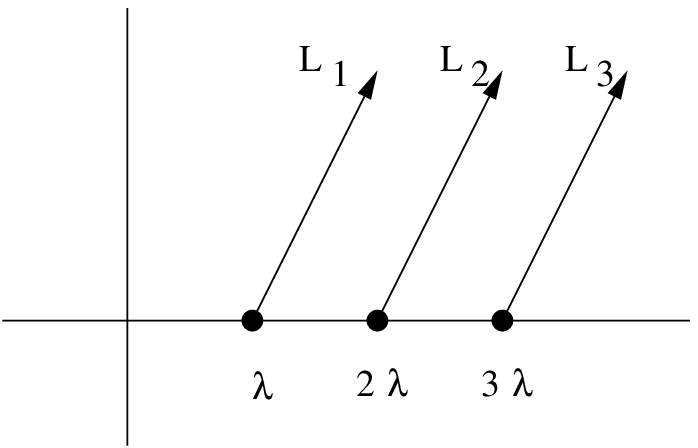}{2in}
$$
The next proposition is a special case of a general result that will
appear in subsequent work of the authors. 

\begin{theorem}
\lbl{thm.transseries}
\rm{(a)}For every $k \in \BN^+$, there exist analytic integable functions 
$R_k \in L^1[0,\infty)$ such that:
\begin{equation}
b_n=\lambda^{-n} n^{3/2} \sum_{k=1}^\infty \frac{1}{k^n} \int_0^\infty
e^{-n p} R_k(e^p) dp
\end{equation}
\rm{(b)}Moreover, for every $k \in \BN^+$, we have an asymptotic expansion
\begin{equation}
\lbl{eq.Gj}
\int_0^\infty e^{-np} R_k(p) dp \sim \sum_{l=0}^\infty \frac{c_{k,l}}{n^l} p^l.
\end{equation}
\rm{(c)} Thus, $G(p)$ determines the transseries \eqref{eq.Gformula4}.
Conversely, $G(p)$ is uniquely determined by its transseries.
\end{theorem}
The functions $R_k$ are constructed from the jump (i.e., variation) of the
multivalued function $G(p)$ at the rays $L_k$.

\begin{proof}
The proof is a well-known application of Cauchy's formula and a deformation
of the contour; 
see for example \cite{Ju}. For the benefit of the reader, we give the details. 
For a technical integrability reason we will work with the following variation
$g(p)$ of $G(p)$:

\begin{equation}
\lbl{eq.gG}
g(p)=\sum_{n=1}^\infty \frac{b_n}{n^2} p^n
\end{equation}
which of course satisfies
$$
\left(p \frac{d}{dp}\right)^2 g(p)=G(p)-b_0.
$$
Of course $g$ and $G$ have the same singularities. 
Since $g(p)$ is analytic in a neighborhood of zero, Cauchy's formula
implies that
$$
\frac{b_n}{n^2}=\frac{1}{2 \pi i} \int_{\gamma} \frac{g(p)}{p^{n}}\frac{dp}{p}.
$$
Now, we will deform the contour $\gamma$ in the following way. Choose 
(Hankel) contours $C_k$ along each ray $L_k$, and choose a truncation
$C_k^r$ of them for $r$ large. Join $\cup_{k=1}^r C_1^r$ together as shown
in Figure \ref{f.deformedcontour} for $r=3$, and create a contour $\gamma_r$
\begin{figure}[htpb]
$$
\psdraw{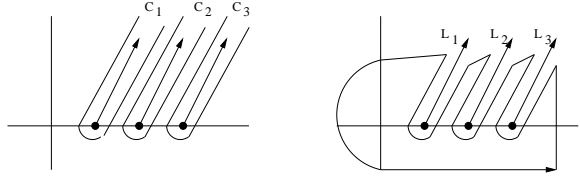}{3in}
$$
\caption{On the left, Hankel contours $C_k$ around each ray $L_k$,
oriented counterclockwise. On the right, a 
truncated contour.}\lbl{f.deformedcontour}
\end{figure}
For every $r$, there is a deformation of $\gamma$ to $\gamma_r$ which does
not pass through the singularities of $g(p)$. It follows that

$$
\frac{b_n}{n^2}
=\frac{1}{2 \pi i} \sum_{k=1}^r \int_{C_k^r} \frac{g(p)}{p^{n}}\frac{dp}{p}+
\frac{1}{2 \pi i} \int_{\Gamma^r} \frac{g(p)}{p^{n}}\frac{dp}{p}
$$
where $\Gamma_r= \gamma_r-\cup_{k=1}^r C_1^r$ is the part of the contour
$\gamma_r$ that is not included in the truncated Hankel contours.
Now, let $r \to \infty$. An estimate shows that 
$$
\lim_{r \to \infty} \int_{\Gamma^r} \frac{g(p)}{p^{n}}\frac{dp}{p}=0.
$$
Let $H_k$ denote a Hankel contour around the ray $L_k$. For every 
$k \in \BN^+$, we have

\begin{equation}
\lbl{eq.CH}
\lim_{r \to \infty} \frac{1}{2 \pi i}
\int_{C_k^r} \frac{g(p)}{p^{n}}\frac{dp}{p}=
\int_{H_k} \frac{g(p)}{p^{n}}\frac{dp}{p}.
\end{equation}
Recall that $g(p)$ is analytic in $\BC\setminus \cup_{k \in \BN^+} L_k$.
For $p \in L_k$, we define the jump (i.e., the variation) $g_k(p)$
of $g(p)$ by

\begin{equation}
\lbl{eq.var}
g_k(x)=\lim_{\e\to 0^+} g(p+i\e)-g(p-i\e).
\end{equation}
On the other hand, Theorem \ref{thm.2} implies that around $p=k\l$, $g$
has an expansion of the form
$$
g(p)=\frac{S_k(p-k\l)}{(p-k\l)^{1/2}}
$$
where $S_k$ is analytic and integrable in $[0,\infty)$. 
It follows that for $p \in L_k$ we have
$$
g_k(p)=2 \frac{S_k(p-k\l)}{(p-k\l)^{1/2}}
$$
Thus, for $t \in \BR^+$ we can write
$$
g_k(k\l e^t)=\frac{T_k(t)}{t^{1/2}}
$$
where $T_k(t)$ is analytic and integrable at $[0,\infty)$. 
A change of variables $p=k\l e^{t+i\theta}$ in Equation \eqref{eq.CH} gives 

$$
\int_{H_k} \frac{g(p)}{p^{n}}\frac{dp}{p}=
(k\l)^{-n} \int_0^\infty e^{-nt}\frac{T_k(t)}{t^{1/2}} dt.
$$
Since

$$
\int_0^\infty e^{-nt} t^c dt=\frac{\Gamma(c+1)}{n^{1+c}}
$$
for all $c \in \BC$ with $\Re(c) \geq -1/2$, it follows that we can write

$$
\int_0^\infty e^{-nt}\frac{T_k(t)}{t^{1/2}} dt=n^{-1/2}
\int_0^\infty e^{-nt} R_k(t) dt
$$
for $R_k$ analytic and integrable at $[0,\infty)$. This proves part (a).

Part (b) follows from Watson's lemma; see \cite{O}.

Part (c) also follows from Watson's lemma.
\end{proof}

\begin{remark}
\lbl{rem.transseries}
As is obvious from the statement and the proof, Theorem 
\ref{thm.transseries} holds for a wide class of resurgent functions $G(p)$,
that includes all \srb\ functions with singularities in a finite set
of rays $ \lambda_1 \BN^+ \cup \dots \cup \lambda_r \BN^+$.
\end{remark}

Among other things, the above theorem makes clear the usefulness (and the
necessity) of transseries versus asymptotic expansions. The asymptotic
expansion \eqref{eq.Gformula3} determines $G(p)$ modulo {\em exponentially
small corrections}. These corrections, beyond all orders in $1/n$, are
precisely captured by the transseries. Theorem \ref{thm.transseries} gives
a synthesis of $G(p)$ by its transseries.
In addition, Theorem \ref{thm.transseries} gives a proof of Theorem
\ref{thm.11}.

\begin{proof}(of Theorem \ref{thm.11})
Part (a) is a general statement about Laplace transforms, and follows from
Watson's lemma.

Part (b) follows from Theorem \ref{thm.transseries} above, and from
Equation \eqref{eq.conj1}. 

For part (c), $F_{\calK}(x)$ determines (via object synthesis), the
analytic function $S^{\med}_{\calK}(x)$, and its radial limits via
\eqref{eq.conj1}. Conversely, the sequence $(\phi_{\calK}(e^{2 \pi i/n}))$
determines its transseries, which in turn determines (via Theorem 
\ref{thm.transseries}) the function $S^{\med}_{\calK}(x)$, which finally
determines $F_{\calK}(x)$ by \eqref{eq.conj2}.
This completes the proof of Theorem \ref{thm.11}.
\end{proof}

\appendix

\section{Resurgence of the power series of the Poincar\'e homology sphere}
\lbl{sec.appendix}

In this section, let $M$ denote the {\em Poincar\'e homology sphere}, a 
closed 3-manifold. In \cite{LZ}, Lawrence-Zagier compute that 
\begin{equation}
\lbl{eq.S1}
F_M(x)=\sum_{n=0}^\infty \frac{a_n}{n!} \frac{1}{(120 x)^n}
\end{equation}
where
\begin{equation}
\lbl{eq.S2}
\sum_{n=0}^\infty \frac{a_{n}}{(2n)!} p^{2n}=\frac{\cos 5p \cos 9p}{\cos 15 p}.
\end{equation}
A computation analogous to the one in Section \ref{sec.boreltransform} 
shows that the Borel transform $G_M(p)$ of $F_M(x)$ is given by:
\begin{equation}
G_M(p)=
c_1 \sum_{n=0}^\infty \frac{\chi_1(n)}{(-30p+n^2 \pi^2)^{3/2}} +
c_2 \sum_{n=0}^\infty \frac{\chi_2(n)}{(-30p+n^2 \pi^2)^{3/2}}
\end{equation}
where
\begin{equation}
\displaystyle 
c_1 =  \frac{\sqrt{6(5+\sqrt{5})}}{120},
\qquad
c_2 =  \frac{\sqrt{6(5-\sqrt{5})}}{120},
\end{equation}
and $\chi_1$, $\chi_1$ are periodic functions defined by the table:
$$
\begin{tabular}{|l|l|l|l|l|l|l|l|l|l|} \hline
$ n \bmod 60 $ & $7$ & $13$ & $17$ & $23$ 
& $37$ & $43$ & $47$ & $53$ & other \\ \hline 
$ \chi_1(n) $ & $-1$ & $-1$ & $-1$ & $-1$ & $1$ & $1$ & $1$ & $1$ & $0$ 
\\ \hline  
\end{tabular}
$$
and
$$
\begin{tabular}{|l|l|l|l|l|l|l|l|l|l|} \hline
$ n \bmod 60 $ & $1$ & $11$ & $19$ & $29$ 
& $31$ & $41$ & $49$ & $59$ & other \\ \hline 
$ \chi_2(n) $ & $-1$ & $-1$ & $-1$ & $-1$ & $1$ & $1$ & $1$ & $1$ & $0$ 
\\ \hline  
\end{tabular}
$$


\ifx\undefined\bysame
        \newcommand{\bysame}{\leavevmode\hbox
to3em{\hrulefill}\,}
\fi


\begin{thebibliography}{[EMSS]}

\bibitem[C]{C} O. Costin,
        {\em On Borel summation and Stokes phenomena for rank-$1$ nonlinear 
        systems of ordinary differential equations},
        Duke Math. J.  {\bf 93}  (1998) 289--344.

\bibitem[CK]{CK} \bysame and M.D. Kruskal,
        {\em On optimal truncation of divergent series solutions of 
        nonlinear differential systems},
        R. Soc. Lond. Proc. Ser. A Math. Phys. Eng. Sci.  {\bf 455}  (1999)
        no. 1985, 1931--1956.

\bibitem[CG1]{CG1} \bysame and S. Garoufalidis,
        {\em Resurgence of the Euler-MacLaurin summation formula}, 
        Annales de l' Institut Fourier, {\bf 58} (2008) 893--914. 

\bibitem[CG2]{CG2} \bysame and \bysame,
        {\em Resurgence of 1-dimensional series of sum-product type},
        in preparation.

\bibitem[CNP1]{CNP1} B. Candelpergher, J.C. Nosmas and F. Pham, 
        {\em Approche de la r\'esurgence},
        Actualit\'es Math\'ematiques, Hermann 1993.

\bibitem[CNP2]{CNP2} \bysame, \bysame and \bysame,
        {\em Premiers pas en calcul \'etranger},
        Ann. Inst. Fourier (Grenoble) {\bf 43}  (1993) 201--224.
 
\bibitem[D]{D} E. Delabaere,
        {\em Introduction to the \'Ecalle theory},
        in Computer algebra and differential equations,
        London Math. Soc. Lecture Note Ser., {\bf 193} (1994) 59--101.

\bibitem[DP]{DP} \bysame and F. Pham,
        {\em Resurgent methods in semi-classical asymptotics},
        Ann. Inst. H. Poincar\'e Phys. Th\'eor.  {\bf 71}  (1999) 1--94. 

\bibitem[Ec]{Ec} J. \'Ecalle, 
        {\em Resurgent functions}, 
        Vol. I-III Mathematical Publications of Orsay {\bf 81-05} 1981, ibid
        {\bf 81-06} 1981, ibid {\bf 85-05} 1985.

\bibitem[EM]{EM} \bysame and F. Menous,
        {\em  Well-behaved convolution averages and the non-accumulation 
        theorem for limit-cycles},
        in ``The Stokes Phenomenon and Hilbert's 16th problem'',
        World Scientific (1996) 71--102.

\bibitem[GL]{GL} S. Garoufalidis and T.T.Q. Le, 
        {\em Gevrey series in quantum topology},
        J. Reine Angew. Math., (2007) 1--27 {\em in press}.

\bibitem[Ga]{Ga} \bysame,
        {\em Chern-Simons theory, analytic continuation and arithmetic},
        preprint 2007 {\tt arXiv:0711.1716}, 
         Acta Math. Vietnam.  {\bf 33}  (2008) 335--362.

\bibitem[GS]{GS} V. Gelfreich and D. Sauzin, 
        {\em Borel summation and splitting of separatrices for the H\'enon 
        map},
        Ann. Inst. Fourier  {\bf 51} (2001) 513--567.




\bibitem[Ha1]{Ha1} K. Habiro,
        {\em On the quantum $\mathfrak{sl}_2$
        invariants of knots and integral homology
        spheres},
        Geom. Topol. Monogr. {\bf 4} (2002) 55--68.

\bibitem[Ha2]{Ha2} \bysame,
        {\em Cyclotomic completions of polynomial rings},
        Publ. Res. Inst. Math. Sci.  {\bf 40}  (2004) 1127--1146. 

\bibitem[Hi1]{Hi1} K. Hikami,
        {\em Quantum invariant for torus link and modular forms},
        Comm. Math. Phys.  {\bf 246}  (2004)  403--426. 

\bibitem[Hi2]{Hi2} \bysame,
        {\em On the quantum invariant for the Brieskorn homology spheres},
        Internat. J. Math.  {\bf 16}  (2005) 661--685.


\bibitem[HL]{HL} V. Huynh and T.T.Q. Le,
        {\em On the Colored Jones Polynomial and the Kashaev invariant},
        Fundam. Prikl. Mat. {\bf 11} (2005) 57--78.

\bibitem[Ju]{Ju} R. Jungen,
        {\em Sur les s\'eries de Taylor n'ayant que des singularit\'es 
        alg\'ebrico-logarithmiques sur leur cercle de convergence},
        Comment. Math. Helv. {\bf 3} (1931) 266--306. 
 
\bibitem[Ka]{Ka} R. Kashaev,
        {\em The hyperbolic volume of knots from the quantum dilogarithm},
        Modern Phys. Lett. A {\bf 39} (1997) 269--275.

\bibitem[LZ]{LZ} R. Lawrence and D. Zagier, 
        {\em Modular forms and quantum invariants of $3$-manifolds},
        in Sir Michael Atiyah: a great mathematician of the twentieth century,
        Asian J. Math. {\bf 3} (1999) 93--107.

\bibitem[LMO]{LMO}  T.T.Q.  Le, J. Murakami, T. Ohtsuki,
        {\em A universal quantum invariant of 3-manifolds},
        Topology {\bf 37} (1998) 539--574. 

\bibitem[Le]{Le} \bysame,
        {\em Integrality and symmetry of quantum link invariants},
        Duke Math. J.  {\bf 102}  (2000) 273--306. 

\bibitem[Lb]{Lb} N.N. Lebedev, 
        {\em Special functions and their applications},
        Dover Publications, Inc. 1972.

\bibitem[Ml]{Ml}  B. Malgrange, 
        {\em Introduction aux travaux de J. \'Ecalle}, 
        Enseign. Math.   {\bf 31}  (1985) 261--282.
 
\bibitem[Ma]{Ma} S. Mandelbrojt, 
        S\'eries lacunaires, Hermann (1936) pp 18.

\bibitem[M]{M} F. Menous, 
        {\em The well-behaved Catalan and Brownian averages and their 
        applications to real resummation},
        Proceedings of the Symposium on Planar Vector Fields (Lleida, 1996).  
        Publ. Mat.  {\bf 41}  (1997) 209--222.

\bibitem[MM]{MM} H. Murakami and J. Murakami,
        {\em The colored Jones polynomials and the simplicial volume of a 
        knot},
        Acta Math.  {\bf 186}  (2001) 85--104.

\bibitem[O]{O} F. Olver, 
        {\em Asymptotics and special functions}, Reprint. 
        AKP Classics. A K Peters, Ltd., Wellesley, MA, 1997.

\bibitem[OSZ]{OSZ} C. Oliv\'e, D. Sauzin and T.M. Seara,
        {\em Two examples of resurgence},
        in  Analyzable functions and applications,  
        Contemp. Math. {\bf 373} Amer. Math. Soc. (2005) 355--371.

\bibitem[Ra]{Ra} J.P. Ramis, 
        {\em S\'eries divergentes et th\'eories asymptotiques}, 
        Bull. Soc. Math. France  {\bf 121}  (1993),  
        Panoramas et Syntheses, suppl.
 
\bibitem[RT]{RT} N. Reshetikhin, ~V. Turaev,
        {\em Ribbon graphs and their invariants derived from quantum groups},
        Commun.  Math.  Phys.  {\bf 127} (1990) 1--26.

\bibitem[Tu1]{Tu1} V. Turaev,
        {\em The Yang-Baxter equation and invariants of links},
        Inventiones Math. {\bf 92} (1988) 527--553.

\bibitem[Tu2]{Tu2} \bysame,
        {\em Quantum invariants of knots and 3-manifolds}, 
        de Gruyter Studies in Mathematics {\bf 18}, 
        Walter de Gruyter, Berlin New York 1994.

\bibitem[Wi]{Wi} E.Witten, 
        {\em Quantum field theory and the Jones polynomial}, 
        Commun. Math. Physics. {\bf 121} (1989) 360--376.
 
\bibitem[Za]{Za} D. Zagier,
        {\em Vassiliev invariants and a strange identity related to the 
        Dedekind eta-function},
        Topology  {\bf 40}  (2001) 945--960. 

\end{thebibliography}
\end{document}